\documentclass[12pt]{article}

\usepackage{graphicx}
\usepackage{latexsym,amssymb}
\usepackage{amsthm}
\usepackage{indentfirst}
\usepackage{amsmath}
\usepackage{color}
\usepackage{xcolor}
\usepackage{fourier}
\usepackage[colorlinks=true,backref=page]{hyperref}
\usepackage{mathrsfs}
\usepackage{enumitem}
\usepackage{float}

\newtheorem{theoremalph}{Theorem}

\newtheorem*{Main Theorem}{Main Theorem}
\newtheorem{Coro}[theoremalph]{Corollary}
\newtheorem{Theorem}{Theorem}[section]
\newtheorem*{Theorem A}{Theorem A}
\newtheorem*{Theorem A'}{Theorem A'}
\newtheorem*{Theorem B'}{Theorem B'}

\newtheorem*{Conjecture}{Conjecture}
\newtheorem{Definition}[Theorem]{Definition}
\newtheorem{Proposition}[Theorem]{Proposition}
\newtheorem{Lemma}[Theorem]{Lemma}
\newtheorem{Example}[Theorem]{Example}

\newtheorem*{Question*}{Question}

\newtheorem{Remark}[Theorem]{Remark}
\newtheorem*{Remark*}{Remark}

\newtheorem{Remark-numbered}[Theorem]{Remark}

\newtheorem{Claim-numbered}[Theorem]{Claim}

\begin{document}

	\title{Unstable Entropy and Volume Growth for Partially Hyperbolic Diffeomorphisms}

\author{Yuntao Zang\\ \small School of Mathematical Sciences, Soochow University, Suzhou, P.R. China\\
	\small \texttt{ytzang@suda.edu.cn}}

\date{}

	\maketitle
	\begin{abstract}

		\medskip
	Let $f$ be a $C^{1}$ diffeomorphism on a compact manifold $M$ admitting a partially hyperbolic splitting $TM=E^{cs}\oplus E^{u}$ where $E^{u}$ is uniformly expanding. We establish a lower bound for the unstable topological entropy (w.r.t. $E^{u}$): $$h^{u}_{\rm top}\left(f\right)\geq\limsup_{n\to+\infty}\frac{1}{n}\log\int_{M}\left|\det \left(Df^{n}_{x}|_{E^{u}}\right)\right|\,dx.$$ 
	Moreover, if there exists a measure of maximal unstable entropy whose Lyapunov exponents along $E^{cs}$ are all non-positive, then equality holds and the $\limsup$ is in fact a limit.
	
	\medskip
	\noindent\textbf{\normalsize MSC Class:} 37D30; 37B40; 37A35
	\end{abstract}

	\section{Introduction}
	\subsection{Main results}
	Let $f$ be a $C^{1}$ diffeomorphism on a compact Riemannian manifold $M$ admitting a partially hyperbolic splitting $TM=E^{cs}\oplus E^{u}$, where $E^{u}$ is uniformly expanding. Following Hu--Hua--Wu \cite{HHW17}, the \emph{unstable topological entropy} $h^{u}_{\rm top}\left(f\right)$ is defined as the supremum, over $x\in M$, of the topological entropy of local unstable manifolds $W^{u}_{\delta}\left(x\right)$. See Definition \ref{definition of unstable entropy} for the precise definition.
	
	We provide a lower bound for the unstable topological entropy in terms of the integrated volume growth along $E^{u}$.
	
	\begin{theoremalph}\label{thm:main}
		Let $f$ be a $C^{1}$ diffeomorphism on a compact manifold $M$ admitting a partially hyperbolic splitting $TM=E^{cs}\oplus E^{u}$, where $E^{u}$ is uniformly expanding. Then $$h^{u}_{\rm top}\left(f\right)\geq\limsup_{n\to+\infty}\frac{1}{n}\log\int_{M}\left|\det\left(Df^{n}_{x}|_{E^{u}}\right)\right|\,dx.$$ 
		Moreover, if there exists a measure of maximal unstable entropy \footnote{An $f$-invariant probability measure $\mu$ is called a \emph{measure of maximal unstable entropy} if $h^{u}_{\mu}\left(f\right)=h^{u}_{\rm top}\left(f\right)$ where $h^{u}_{\mu}\left(f\right)$ denotes the \emph{unstable metric entropy}. It always exists in our setting. For the definition of unstable metric entropy and the existence of measures of maximal unstable entropy, see Hu--Hua--Wu \cite{HHW17}.} whose Lyapunov exponents along $E^{cs}$ are all non-positive, then equality holds and the $\limsup$ is in fact a limit.
	\end{theoremalph}

Note that the integral is taken with respect to the Lebesgue measure \footnote{Throughout, we assume that the Lebesgue measure on $M$ is normalized to be a probability measure.}. The Lyapunov regular set, i.e., the set of points with well-defined Lyapunov exponents, may have zero Lebesgue measure and hence need not contribute to the integral above. In general, the inequality in Theorem \ref{thm:main} can be strict (see Example \ref{example >}). The underlying mechanism is that a measure of maximal unstable entropy may have positive Lyapunov exponents along $E^{cs}$, so that nearby Lebesgue-typical points are repelled from its support and eventually enter some region where $\left|\det\left(Df^{n}_{x}|_{E^{u}}\right)\right|$ grows at a substantially smaller rate.

\vskip 5pt

The order of "$\log$" and "$\int$" in Theorem \ref{thm:main} is essential. Indeed, Jensen's inequality gives $$\int_{M}\log\left|\det\left(Df^{n}_{x}|_{E^{u}}\right)\right|\,dx\leq\log\int_{M}\left|\det\left(Df^{n}_{x}|_{E^{u}}\right)\right|\,dx.$$ Thus, moving the logarithm inside the integral leads to a smaller asymptotic growth rate. Inspired by Zang \cite[Corollary C]{Zan26}, the following corollary shows that equality between these two growth rates has significant dynamical consequences.

	\begin{Coro}\label{log integeral order}
		Let $f$ be a $C^{1}$ Anosov diffeomorphism on a compact manifold $M$ admitting a hyperbolic splitting $TM=E^{s}\oplus E^{u}$. The following two properties are equivalent: 
		
		\begin{enumerate}[label=(\Roman*)]
			\item $$\limsup_{n\to+\infty}\frac{1}{n}\int_{M} \log\left|\det\left(Df_x^n|_{E_x^u}\right)\right|\,dx=\limsup_{n\to+\infty}\frac{1}{n}\log\int_{M} \left|\det\left(Df_x^n|_{E_x^u}\right)\right|\,dx. $$
			\item There exists a sequence $n_k\to+\infty$ such that, for Lebesgue almost every $x\in M$, every weak-$\ast$ accumulation point of
			$$\left\{\frac{1}{n_k}\sum_{i=0}^{n_k-1}\delta_{f^i\left(x\right)}\right\}_{k\ge1}$$
			is an SRB measure \footnote{An invariant probability measure is called an SRB measure here if it satisfies Pesin's entropy formula, that is, its metric entropy equals the integral of the sum of its positive Lyapunov exponents, counted with multiplicities.} which is also a measure of maximal entropy.
		\end{enumerate}
		
	\end{Coro}

\begin{proof}
	\textbf{Part:} (I) $\Longrightarrow$ (II).
		
	Since $f$ is Anosov, we have $h_{\rm top}^u\left(f\right)=h_{\rm top}\left(f\right)$ and all Lyapunov exponents along $E^s$ are negative. Hence Theorem~\ref{thm:main} and (I) give
	\begin{equation}\label{integeral equal to top}
		\limsup_{n\to+\infty}\frac{1}{n}\int_{M} \log\left|\det\left(Df_x^n|_{E_x^u}\right)\right|\,dx=\limsup_{n\to+\infty}\frac1n\log\int_M\left|\det\left(Df_x^n|_{E_x^u}\right)\right|\,dx=h_{\rm top}\left(f\right).
	\end{equation}
	A general measure-theoretic result (see Lemma \ref{log-integral}) yields a sequence $n_{k}\to +\infty$ such that for Lebesgue almost every $x$, $$\lim_{k\to+\infty}\frac{1}{n_{k}}\log\left|\det\left(Df_x^{n_{k}}|_{E_x^u}\right)\right|=h_{\rm top}\left(f\right).$$ By Catsigeras--Cerminara--Enrich \cite[Corollary 1.8]{CCE17}, for Lebesgue almost every $x$, any weak-$\ast$ limit $\mu$ of
	$$\left\{\frac1{n_k}\sum_{i=0}^{n_k-1}\delta_{f^i\left(x\right)}\right\}_{k\geq1}$$ is an SRB measure. Hence
	$$h_\mu\left(f\right)=\int_M\log\left|\det\left(Df_y|_{E_y^u}\right)\right|\,d\mu\left(y\right)=\lim_{k\to+\infty}\frac{1}{n_{k}}\log\left|\det\left(Df_x^{n_{k}}|_{E_x^u}\right)\right|=h_{\rm top}\left(f\right),$$ which shows that $\mu$ is also a measure of maximal entropy. This proves (II).
	
	\textbf{Part:} (II) $\Longrightarrow$ (I).
	
	For Lebesgue almost every $x$, by the SRB property, we have
	$$\lim_{k\to+\infty}\frac1{n_k}\log\left|\det\left(Df_x^{n_k}|_{E_x^u}\right)\right|=h_{\rm top}\left(f\right).$$ Hence
	$$\lim_{k\to+\infty}\frac1{n_k}\int_M\log\left|\det\left(Df_x^{n_k}|_{E_x^u}\right)\right|\,dx=h_{\rm top}\left(f\right).$$ By Theorem~\ref{thm:main} and Jensen's inequality, we then have
	$$h_{\rm top}\left(f\right)\leq \limsup_{n\to+\infty}\frac1n\int_M\log\left|\det\left(Df_x^n|_{E_x^u}\right)\right|\,dx\leq \limsup_{n\to+\infty}\frac{1}{n}\log\int_{M} \left|\det\left(Df_x^n|_{E_x^u}\right)\right|\,dx\leq h_{\rm top}\left(f\right),$$ which proves (I).
\end{proof}

	Property (II) in Corollary \ref{log integeral order} may appear somewhat mysterious. This reflects a subtle difference between the $C^{1}$ and $C^{2}$ settings for Anosov diffeomorphisms (see Catsigeras-Cerminara-Enrich \cite{CCE17}). For $C^{2}$ Anosov diffeomorphisms, Property (II) can be replaced by the following more natural condition:
	
	\begin{enumerate}[label=(\Roman*)',start=2]
		\item There exist finitely many ergodic SRB measures, each of which is a measure of maximal entropy, and the union of their basins has full Lebesgue measure.
	\end{enumerate}
	
	\vskip 5pt
	
\subsection{Applications}

We apply Theorem~\ref{thm:main} to the following four classes of partially hyperbolic diffeomorphisms. In cases (2)--(4), the splitting is $TM=E^s\oplus E^c\oplus E^u$, and we set $E^{cs}:=E^s\oplus E^c$.

\begin{enumerate}[label=(\arabic*)]
	\item \emph{Maps factoring over Anosov with $c$-mostly contracting center.}
	Consider $C^1$ dynamically coherent partially hyperbolic diffeomorphisms with splitting $TM=E^{cs}\oplus E^u$ which factor over a hyperbolic linear toral automorphism $A$ and have $c$-mostly contracting center, in the sense of \cite{UVYY24}.
	
	In particular, this class includes \emph{derived-from-Anosov diffeomorphisms with contracting linear center} (see \cite[Section 6]{UVYY24} for the precise setting and its properties).
	
	\item \emph{Partially hyperbolic diffeomorphisms of three-dimensional nilmanifolds.}
	Consider $C^2$ partially hyperbolic diffeomorphisms of a compact three-dimensional nilmanifold not homeomorphic to $\mathbb T^3$, with one-dimensional stable, center, and unstable bundles; see \cite{UVY21}.
	
	\item \emph{Accessible systems with compact center leaves.}
	Consider $C^{1+\alpha}$ partially hyperbolic diffeomorphisms of a closed three-dimensional manifold which are accessible and dynamically coherent, with compact one-dimensional center leaves; see \cite{RHRHTU12}.
	
	\item \emph{Time-one maps of Anosov flows and certain nearby systems.}
	Let $\varphi^t$ be a $C^2$ transitive Anosov flow on a compact manifold $M$. Let $\mathcal U\subset\operatorname{Diff}^1(M)$ be the open set provided by \cite[Theorem 1.1]{BFT22}, whose closure contains $\varphi^1$. We consider $\varphi^1$ and the partially hyperbolic diffeomorphisms $f\in\mathcal U\cap\operatorname{Diff}^2(M)$, with their one-dimensional center bundles.
\end{enumerate}

\begin{Coro}\label{application}
	For every diffeomorphism $f$ in any of the four classes above, we have
	$$h_{\rm top}^u\left(f\right)=\lim_{n\to+\infty}\frac{1}{n}\log\int_{M}\left|\det\left(Df^{n}_{x}|_{E^{u}}\right)\right|\,dx.$$
	In cases (2)--(4), this common value also equals $h_{\rm top}\left(f\right)$. In case (1), it equals $h_{\rm top}\left(A\right)$.
\end{Coro}

\begin{proof}
	We first note that an ergodic measure of maximal entropy with non-positive exponents along $E^{cs}$ also maximizes unstable entropy \cite[Corollary A.1]{HHW17}.
	
	We now verify the equality criterion in Theorem~\ref{thm:main} for each class.
	\begin{enumerate}[label=(\arabic*)]
		\item By \cite[Theorem B and Proposition 4.2]{UVYY24}, every ergodic measure of maximal unstable entropy has only negative Lyapunov exponents along $E^{cs}$. 
		
		\item By \cite[Theorem A]{UVY21}, there exists an ergodic measure of maximal entropy with non-positive exponent along $E^{cs}$.
		
		\item By \cite[Theorem 1]{RHRHTU12}, there exists an ergodic measure of maximal entropy with non-positive exponent along $E^{cs}$.
		
		\item By \cite[Theorem 1.1]{BFT22}, there exists an ergodic measure of maximal entropy with non-positive exponent along $E^{cs}$ (this also holds for the time-one map).
	\end{enumerate}
	
	In all four cases there exists an ergodic measure of maximal unstable entropy whose Lyapunov exponents along $E^{cs}$ are non-positive. Theorem~\ref{thm:main} therefore gives the asserted integral formula. The preceding observations also establish $h_{\rm top}^u\left(f\right)=h_{\rm top}\left(f\right)$ in cases (2)--(4). Finally, \cite[Lemma 3.1]{UVYY24} gives $h_{\rm top}^u\left(f\right)=h_{\rm top}\left(A\right)$ in case (1).
\end{proof}

	\subsection{Background}
	
	The relationship between entropy and the integrated volume growth has a long history, dating back to Sacksteder--Shub \cite{SS78}. For $C^{r}$ ($r>1$) diffeomorphisms, Przytycki \cite{Prz80} established an upper bound for the topological entropy in terms of the integrated volume growth of the exterior powers of the derivative. Newhouse \cite{New88} extended related estimates to smooth maps. On the other hand, Yomdin \cite{Yom87} developed a fundamental theory relating entropy to the volume growth of smooth submanifolds. Building on these ideas, Kozlovski \cite{Koz98} (see also Klapper--Young \cite{KY95}) established an integral formula for the topological entropy of $C^{\infty}$ maps. 
	
	The integrated volume growth for uniformly hyperbolic systems were considered by Kifer--Newhouse \cite{KN91}. In the partially hyperbolic setting, several results have extended these connections to systems of lower regularity.  Saghin \cite{Sag14} obtained an upper bound for the metric entropy of $C^{1}$ partially hyperbolic diffeomorphisms by combining the integrated volume growth along the unstable bundle with the maximal Lyapunov exponent in the center-stable direction. Guo--Liao--Sun--Yang \cite{GLSY18} established a similar result under dominated splittings. Subsequently, Yang--Zang \cite{YZ23} proved an integral formula for the topological entropy of $C^{1}$ partially hyperbolic diffeomorphisms whose intermediate bundles are one-dimensional. More recently, Zang \cite{Zan26} obtained such a formula only involving the derivative itself for $C^{r}$ ($r>1$) surface diffeomorphisms with sufficiently large topological entropy. Related developments in random dynamical systems can be found in Ma \cite{Ma20} and Ma--Zhang \cite{MZ23}.
	
	The volume growth of unstable manifolds provides another natural approach to measuring dynamical complexity. Cogswell \cite{Cog00} related the metric entropy of an ergodic measure to the volume growth of its unstable manifolds; see also Zang \cite{Zan22}. For partially hyperbolic diffeomorphisms, Hu--Hua--Wu \cite{HHW17} introduced unstable topological entropy, established its variational principle, and proved that it coincides with the unstable volume growth introduced by Hua--Saghin--Xia \cite{HSX08}.
	
	Unlike the preceding integral estimates for the usual topological or metric entropy, our result concerns the \emph{unstable topological entropy} and involves only the Jacobian along the uniformly expanding bundle $E^{u}$. In particular, the inequality in Theorem \ref{thm:main} holds in the $C^{1}$ setting without any restriction on the Lyapunov exponents along $E^{cs}$. Moreover, combining our result with the characterization of unstable topological entropy by Hu--Hua--Wu \cite{HHW17}, we obtain $$\lim_{n\to+\infty}\frac{1}{n}\log\int_{M}\left|\det\left(Df^{n}_{x}|_{E^{u}}\right)\right|\,dx\leq h^{u}_{\rm top}\left(f\right)=\chi_{u}\left(f\right),$$ where the unstable volume growth $\chi_{u}\left(f\right)$ is defined, for any fixed $\delta>0$, by $$\chi_{u}\left(f\right):=\sup_{x\in M}\limsup_{n\to+\infty}\frac{1}{n}\log {\rm Vol}\left(f^{n}\left(W^{u}_{\delta}\left(x\right)\right)\right),$$ and is independent of the choice of $\delta$. Saghin \cite{Sag14} pointed out that the above inequality can be strict (see Example \ref{example >}). On the other hand, he proved that the two growth rates coincide if there exists a closed $u$-form, where $u=\dim E^{u}$, which is non-degenerate on the unstable foliation. Our equality criterion provides another sufficient condition for the two growth rates to coincide. Notice that our condition does not imply $h^{u}_{\rm top}\left(f\right)=h_{\rm top}\left(f\right)$: a measure of maximal unstable entropy need not be a measure of maximal entropy, and additional entropy may arise from invariant measures with positive Lyapunov exponents along $E^{cs}$.
	
	The main new ingredient in the proof of the inequality in Theorem \ref{thm:main} is the following uniformity result: for any sufficiently small $\delta>0$, $$\chi_{u}\left(f\right)=\sup_{x\in M}\limsup_{n\to+\infty}\frac{1}{n}\log{\rm Vol}\left(f^{n}\left(W^{u}_{\delta}\left(x\right)\right)\right)=\lim_{n\to+\infty}\frac{1}{n}\log\sup_{x\in M}{\rm Vol}\left(f^{n}\left(W^{u}_{\delta}\left(x\right)\right)\right).$$ Saghin \cite{Sag14} obtained this equality under the existence of a closed $u$-form non-degenerate on the unstable foliation, while Proposition \ref{sup and lim order} requires no such assumption. The inequality in Theorem \ref{thm:main} then follows from a classical Fubini argument.
	
\subsection{Further questions}	
	It remains unclear whether the $\limsup$ in Theorem \ref{thm:main} is always a limit.
	
	\begin{Question*}
		Let $f$ be a $C^{1}$ diffeomorphism on a compact manifold $M$ admitting a partially hyperbolic splitting $TM=E^{cs}\oplus E^{u}$, where $E^{u}$ is uniformly expanding.
		
		Does the limit $$\lim_{n\to+\infty}\frac{1}{n}\log\int_{M}\left|\det\left(Df^{n}_{x}|_{E^{u}}\right)\right|\,dx$$ always exist?
	\end{Question*}
	
	\begin{Remark*}
		The limit exists if there is a measure of maximal unstable entropy whose Lyapunov exponents along $E^{cs}$ are non-positive; see the proof of Theorem \ref{thm:main}. In particular, this holds for uniformly hyperbolic systems. We do not know whether the limit exists in general.
	\end{Remark*}
	\vskip 5pt
	\begin{Conjecture}
		Let $f$ be a $C^{1}$ diffeomorphism on a compact manifold $M$ admitting a partially hyperbolic splitting $TM=E^{s}\oplus E^{c}\oplus E^{u}$, where $E^{u/s}$ is uniformly expanding / contracting and $\dim E^{c}=1$. Then $$h^{u}_{\rm top}\left(f\right)=\limsup_{n\to+\infty}\frac{1}{n}\log\int_{M}\left|\det\left(Df^{n}_{x}|_{E^{u}}\right)\right|\,dx.$$
	\end{Conjecture}
	\begin{Remark*}
		Example \ref{example >} has center dimension two and it seems that this construct prevents center dimension one. The examples in Corollary \ref{application} also provide further evidence for this conjecture.
	\end{Remark*}

\section{Preliminaries}

\begin{Definition}
	Let $f$ be a $C^{1}$ diffeomorphism on a compact Riemannian manifold $M$. We say that $f$ is \emph{partially hyperbolic} if there exists a $Df$-invariant splitting $$TM=E^{cs}\oplus E^{u}$$ and constants $C>0$ and $0<\lambda<1$ such that, for every $x\in M$ and $n\geq1$, $$\frac{\left\|Df^{n}_{x}|_{E^{cs}}\right\|}{m\left(Df^{n}_{x}|_{E^{u}}\right)}\leq C\lambda^{n},\qquad m\left(Df^{n}_{x}|_{E^{u}}\right)\geq C^{-1}\lambda^{-n},$$ where $m\left(A\right):=\left\|A^{-1}\right\|^{-1}$ denotes the conorm.
\end{Definition}

The bundle $E^{u}$ is called \emph{uniformly expanding} and uniquely integrates to the unstable foliation $W^{u}$. For $x\in M$ and $\delta>0$, let $W^{u}_{\delta}\left(x\right)$ be the $\delta$-ball in $W^{u}\left(x\right)$ with respect to the intrinsic Riemannian distance $d^{u}$. Define $$d^{u}_{n}\left(y,z\right):=\max_{0\leq j<n}d^{u}\left(f^{j}\left(y\right),f^{j}\left(z\right)\right).$$ Let $s^{u}\left(n,\varepsilon,x,\delta\right)$ denote the maximal cardinality of an $\left(n,\varepsilon\right)$ $u$-separated subset of $W^{u}_{\delta}\left(x\right)$.

\begin{Definition}[\cite{HHW17}]\label{definition of unstable entropy}
	The \emph{unstable topological entropy} of $f$ is defined by $$h^{u}_{\rm top}\left(f\right):=\lim_{\delta\to0}\sup_{x\in M}h_{\rm top}\left(f,W^{u}_{\delta}\left(x\right)\right),$$ where $$h_{\rm top}\left(f,W^{u}_{\delta}\left(x\right)\right):=\lim_{\varepsilon\to0}\limsup_{n\to+\infty}\frac{1}{n}\log s^{u}\left(n,\varepsilon,x,\delta\right).$$
\end{Definition}

Equivalently, one may use $\left(n,\varepsilon\right)$ $u$-spanning sets or open covers in the above definition, and the resulting quantity is independent of the choice of $\delta$ (Hu--Hua--Wu \cite{HHW17}).

Let $\operatorname{Vol}$ denote the induced Riemannian volume on submanifolds. The unstable topological entropy admits the following characterization in terms of unstable volume growth.

\begin{Proposition}[Hu--Hua--Wu \cite{HHW17}]\label{unstable entropy equal volume growth}
	For any sufficiently small $\delta>0$, $$h^{u}_{\rm top}\left(f\right)=\sup_{x\in M}\limsup_{n\to+\infty}\frac{1}{n}\log \operatorname{Vol}\left(f^{n}\left(W^{u}_{\delta}\left(x\right)\right)\right).$$
\end{Proposition}
	
\section{Volume growth of the local unstable manifolds}
\subsection{Uniform volume growth}
A key ingredient in the proof of Theorem \ref{thm:main} is the following uniformity property for unstable volume growth.
\begin{Proposition}\label{sup and lim order}
	Let $f$ be a $C^{1}$ diffeomorphism on a compact manifold $M$ admitting a partially hyperbolic splitting $TM=E^{cs}\oplus E^{u}$ where $E^{u}$ is uniformly expanding. Then for any $\delta>0$, $$\sup_{x}\limsup_{n\to+\infty}\frac{1}{n}\log{\rm Vol}\left(f^{n}\left(W^{u}_{\delta}\left(x\right)\right)\right)=\lim_{n\to+\infty}\frac{1}{n}\log \sup_{x}{\rm Vol}\left(f^{n}\left(W^{u}_{\delta}\left(x\right)\right)\right).$$	
\end{Proposition}
\begin{proof}
	We only need to prove the result for all sufficiently small $\delta>0$. Put
	$$h:=\sup_{x\in M}\limsup_{n\to+\infty}\frac{1}{n}\log\operatorname{Vol}\left(f^n\left(W_\delta^u\left(x\right)\right)\right),\qquad A_n:=\sup_{x\in M}\operatorname{Vol}\left(f^n\left(W_\delta^u\left(x\right)\right)\right).$$
	
	We first fix several uniform constants arising from the uniform local geometry of the local unstable manifolds. For simplicity, we refer to a local unstable manifold as a \emph{plaque}.
	
	\begin{itemize}
		\item  Fix an integer $K\geq 1$ such that every $W_{3\delta}^u\left(x\right)$ can be covered by at most $K$ plaques of radius $\delta$.
		\item Fix some constant $C\ge1$ such that, for every $x\in M$ and $m\ge 0$, the set $f^m\left(W_\delta^u\left(x\right)\right)$ can be covered by at most $C\operatorname{Vol}\left(f^m\left(W_{2\delta}^u\left(x\right)\right)\right)$ plaques of radius $\delta$.
	\end{itemize}
	Since $W_{2\delta}^u\left(x\right)$ can be covered by at most $K$ plaques of radius $\delta$, we have
	$$\operatorname{Vol}\left(f^m\left(W_{2\delta}^u\left(x\right)\right)\right)\le K A_m.$$
	Thus $f^m\left(W_\delta^u\left(x\right)\right)$ can be covered by plaques $W_\delta^u\left(y_1\right),\ldots,W_\delta^u\left(y_q\right)$, where $q\leq CK A_m$. Applying $f^n$ to this cover gives
	$$\operatorname{Vol}\left(f^{m+n}\left(W_\delta^u\left(x\right)\right)\right)\leq\sum_{j=1}^q\operatorname{Vol}\left(f^n\left(W_\delta^u\left(y_j\right)\right)\right)\leq qA_n\leq CK A_mA_n.$$
	Taking the supremum over $x$, we obtain
	$$A_{m+n}\leq CK A_mA_n.$$
	Hence $\left\{\log\left(CKA_n\right)\right\}_{n\ge1}$ is subadditive, so the limit
	$$\lim_{n\to+\infty}\frac{1}{n}\log A_n$$
	exists. By definition, this limit is at least $h$.
	
	We prove the reverse inequality using a covering tree. Since $W_{3\delta}^u\left(x\right)$ can be covered by at most $K$ plaques of radius $\delta$, we have
	$$\limsup_{n\to+\infty}\frac{1}{n}\log\operatorname{Vol}\left(f^n\left(W_{3\delta}^u\left(x\right)\right)\right)\le h.$$
	Fix $s>h$, noting that $h\ge0$ by uniform expansion. For each $x\in M$, choose sufficiently large $n_x \in\mathbb{N}$ such that
	$$\operatorname{Vol}\left(f^{n_x}\left(W_{3\delta}^u\left(x\right)\right)\right)<e^{\frac{s+h}{2} n_x},\quad C\cdot e^{\frac{s+h}{2}n_x}\leq e^{sn_x}.$$
	The continuous dependence of local unstable plaques in the $C^1$ topology, together with the radius margin $2\delta<3\delta$, gives a neighborhood $U_x$ of $x$ such that
	$$\operatorname{Vol}\left(f^{n_x}\left(W_{2\delta}^u\left(y\right)\right)\right)<e^{\frac{s+h}{2} n_x}\qquad\left(y\in U_x\right).$$
	Choose a finite subcover $U_{x_1},\ldots,U_{x_\ell}$ and set $N:=\max_i n_{x_i}$. Thus, whenever $y\in U_{x_i}$, the set $f^{n_{x_i}}\left(W_\delta^u\left(y\right)\right)$ admits a cover by at most $C e^{\frac{s+h}{2} n_{x_{i}}}\leq e^{s n_{x_i}}$ plaques of radius $\delta$.
	
	Fix $x\in M$ and $n\ge1$. We construct a rooted tree whose vertices are labelled by plaques $W_\delta^u\left(y\right)$ and carry accumulated times $t_{y}$. The root is labelled by $W_\delta^u\left(x\right)$ and has time zero $t_{x}=0$.
	
	At a vertex $W_\delta^u\left(y\right)$, stop if $n-t_y<N$. Otherwise, choose one index $i$ such that $y\in U_{x_i}$ and fix a cover
	$$f^{n_{x_i}}\left(W_\delta^u\left(y\right)\right)\subset\bigcup_{j=1}^{q}W_\delta^u\left(z_j\right),\qquad q\le e^{s n_{x_i}}.$$
	Attach $q$ new vertices $W_\delta^u\left(z_{1}\right),\ldots,W_\delta^u\left(z_{q}\right)$ as the children of $W_\delta^u\left(y\right)$, setting
	$$t_{z_j}:=t_y+n_{x_i}\qquad\left(1\le j\le q\right).$$ We assume that each occurrence of a plaque is treated as a distinct vertex, with its own accumulated time. Repeat this construction at every non-terminal vertex. Note that the tree is finite, and every terminal vertex $W_\delta^u\left(z\right)$ satisfies $0\le n-t_z<N$.
	
	By a simple induction, the tree has at most $e^{sn}$ terminal vertices. The covering relation implies that $f^n\left(W_{\delta}^u\left(x\right)\right)$ can be covered by at most $e^{sn}\cdot C_{N}$ plaques of radius $\delta$, where each terminal plaque is iterated for fewer than $N$ additional steps, and the additional covering number is bounded by some constant $C_{N}$. By the arbitrariness of $x$,
	$$A_n\le e^{sn}\cdot C_{N}.$$
	Letting $n\to+\infty$ and then $s\downarrow h$, we obtain
	$$\lim_{n\to+\infty}\frac{1}{n}\log A_n\le h,$$ which completes the proof.
\end{proof}
\begin{Remark}\noindent

	\begin{enumerate}
		\item The same argument also yields the corresponding \emph{$\liminf$} version of Proposition \ref{sup and lim order}: for any $\delta>0$, $$\sup_{x}\liminf_{n\to+\infty}\frac{1}{n}\log{\rm Vol}\left(f^{n}\left(W^{u}_{\delta}\left(x\right)\right)\right)=\lim_{n\to+\infty}\frac{1}{n}\log \sup_{x}{\rm Vol}\left(f^{n}\left(W^{u}_{\delta}\left(x\right)\right)\right).$$	
		\item In the proof of Proposition~\ref{sup and lim order}, uniform expansion is used to obtain the required plaque-covering estimates. It may be possible to establish analogous estimates in certain dominated settings by considering a suitable family of submanifolds tangent to cones around the dominating bundle. 
		\item It is also interesting to compare Proposition~\ref{sup and lim order} with the classical formula (see, e.g., Schreiber \cite{Sch98}) for continuous subadditive functions $\left\{\varphi_n\right\}$:
		$$\sup_{x\in M}\limsup_{n\to+\infty}\frac1n\varphi_n\left(x\right)=\lim_{n\to+\infty}\frac1n\sup_{x\in M}\varphi_n\left(x\right).$$ Our initial (failed) attempt was to deduce Proposition~\ref{sup and lim order} from this formula by finding a suitable subadditive formulation of the plaque-volume growth.
	\end{enumerate}

\end{Remark}

\vskip 5pt

For $\eta>0$, define the \emph{unstable cone at $x$} by
$$\mathcal C_\eta^u\left(x\right):=\left\{v_u+v_{cs}:v_u\in E_x^u,\ v_{cs}\in E_x^{cs},\ \left\|v_{cs}\right\|\leq\eta\left\|v_u\right\|\right\}.$$
A submanifold $D$ is \emph{tangent to $\mathcal C_\eta^u$} if $T_xD\subset\mathcal C_\eta^u\left(x\right)$ for every $x\in D$.

Let $B^{u}_r\left(0\right)\subset E_x^u$ denote the ball at $0$ with radius $r$ in $E_x^u$. For sufficiently small $r,\eta>0$, let $\mathcal D_r\left(\eta\right)$ consist of the disks
$$D=\exp_x\left\{v+g\left(v\right):v\in B^{u}_r\left(0\right)\subset E_x^u\right\},$$
where $x\in M$, the map $g: E_x^u\to E_x^{cs}$ is $C^1$, $g\left(0\right)=0$ and the exponential image of the graph of $g|_{B_{10r}^u\left(0\right)}$ is tangent to $\mathcal C_\eta^u$.

We next show that the unstable volume growth can equivalently be computed using disks of a fixed size tangent to a sufficiently narrow unstable cone field.
\begin{Proposition}\label{thm:cone-volume-growth}
	Let $f$ be a $C^1$ diffeomorphism of a compact Riemannian manifold $M$ admitting a partially hyperbolic splitting $TM=E^{cs}\oplus E^u$, where $E^u$ is uniformly expanding. For all small $\delta,\eta,r>0$,
	$$\lim_{n\to+\infty}\frac{1}{n}\log\sup_{x\in M}\operatorname{Vol}\left(f^n\left(W_\delta^u\left(x\right)\right)\right)=\lim_{n\to+\infty}\frac{1}{n}\log\sup_{D\in\mathcal D_r\left(\eta\right)}\operatorname{Vol}\left(f^{n}\left(D\right)\right).$$
\end{Proposition}

\begin{proof}
	This is essentially a consequence of domination. Fix small $r>0$ so that the exponential map is nearly isometric. Write $$h:=\lim_{m\to+\infty}\frac{1}{m}\log A_{m}\quad \text{ where } \quad A_m:=\sup_{x\in M}\operatorname{Vol}\bigl(f^m\left(W_{5r}^u\left(x\right)\right)\bigr).$$  Note that these graphs in $\mathcal D_r\left(\eta\right)$ can be arbitrarily $C^1$-close to the local unstable manifolds, provided $\eta>0$ is sufficiently small. Consequently, given any large $m$, we can choose $\eta_{m}>0$ sufficiently small such that 
	$$\sup_{D\in\mathcal D_{r}\left(\eta_{m}\right)}\operatorname{Vol}\left(f^m \left(\widehat D\right)\right)\leq 2A_m,$$ where $\widehat D$ is the $2r$-extension of $D$:
	$$\widehat D:=\exp_x\left\{v+g\left(v\right):v\in B_{2r}^u\left(0\right)\right\}.$$
	By domination, there are $C\geq 1$ independent of $m$ and $C_{m}\geq 1$ such that
	
	\begin{enumerate}
		\item for any $D\in\mathcal D_r\left(\eta_{m}\right)$, $$
		f^m\left(D\right)\subset\bigcup_{j=1}^qD_j,\qquad D_j\in\mathcal D_r\left(\eta_m\right),\qquad q\leq C\operatorname{Vol}\left(f^m\widehat D\right)\leq 2CA_{m},$$
		\item for any $D\in\mathcal D_r\left(\eta_{m}\right)$ and any $0\leq j\leq m-1$, $\operatorname{Vol}\left(f^{j}\left(D\right)\right)\leq C_m$,
		\item By domination, there is some large $k_{m}\geq 1$ such that $$\sup_{D\in\mathcal D_r\left(\eta\right)} \operatorname{Vol}\left(f^{n+k_m}\left(D\right)\right) \leq C_m\sup_{D\in\mathcal D_r\left(\eta_m\right)} \operatorname{Vol}\left(f^n\left(D\right)\right). $$
	\end{enumerate}
	For $n\geq k_m$, write $n=k_m+qm+j$, where $0\leq j<m$. Repeating the $m$-step covering and using the finite-time bound for the remaining $j$ iterates, we obtain
	$$\sup_{D\in\mathcal D_r\left(\eta\right)}\operatorname{Vol}\left(f^n\left(D\right)\right)\leq C_m^2\left(2CA_m\right)^q.$$
	Since $q/n\to1/m$ as $n\to+\infty$, it follows that
	$$\limsup_{n\to+\infty}\frac1n\log\sup_{D\in\mathcal D_r\left(\eta\right)}\operatorname{Vol}\left(f^n\left(D\right)\right)\leq\frac1m\log\left(2CA_m\right).$$

Letting $m\to+\infty$ gives the upper bound $h$. The reverse inequality, $$\liminf_{n\to+\infty}\frac{1}{n}\log\sup_{D\in\mathcal D_r\left(\eta\right)}\operatorname{Vol}\left(f^n\left(D\right)\right)\geq\liminf_{n\to+\infty}\frac{1}{n}\log\sup_{x\in M}\operatorname{Vol}\left(f^n\left(W_\delta^u\left(x\right)\right)\right)=h$$ is an immediate consequence of the uniform finite covering of unstable plaques, which completes the proof.
\end{proof}

\subsection{Proof of Theorem \ref{thm:main}}
The following estimate is essentially a Fubini-type argument. However, in the $C^{1}$ setting, the unstable foliation need not be absolutely continuous, so one cannot directly disintegrate the ambient volume along unstable leaves. We therefore work instead with disks tangent to a narrow unstable cone field and use Proposition \ref{thm:cone-volume-growth} to recover the same unstable volume growth rate.
\begin{Lemma}\label{unstable entropy geq integrated volume growth}
	Let $f$ be a $C^{1}$ diffeomorphism on a compact manifold $M$ admitting a partially hyperbolic splitting $TM=E^{cs}\oplus E^{u}$ where $E^{u}$ is uniformly expanding. For any $\delta>0$, $$\lim_{n\to+\infty}\frac{1}{n}\log \sup_{x}{\rm Vol}\left(f^{n}\left(W^{u}_{\delta}\left(x\right)\right)\right)\geq\limsup_{n\to+\infty}\frac{1}{n}\log\int_{M}\left|\det \left(Df^{n}_{x}|_{E^{u}}\right)\right|\,dx.$$	
\end{Lemma}

\begin{proof}
	Fix any small $r$ and $\eta$. Write
	$$J_n^u\left(x\right):=\left|\det\left(Df_x^n|_{E_x^u}\right)\right|,\qquad A_n:=\sup_{D\in\mathcal D_r\left(\eta\right)}\operatorname{Vol}\left(f^{n}\left(D\right)\right).$$
	
	We choose finitely many sufficiently small $C^1$ product charts $\left\{\Phi_i:E^{u}_{x_{i}}\times E^{cs}_{x_{i}}\to M\right\}_{1\leq i\leq N}$ covering $M$ and whose horizontal slices are contained in disks from $\mathcal D_r\left(\eta\right)$. More precisely, writing
	$$S_{i,w}:=\Phi_i\left(B_{r/10}^u\left(0\right)\times\left\{w\right\}\right),$$
	for small $w$, there exists $D_{i,w}\in\mathcal D_r\left(\eta\right)$ such that $S_{i,w}\subset D_{i,w}$. Consequently,
	$$\operatorname{Vol}\left(f^n\left(S_{i,w}\right)\right)\leq\operatorname{Vol}\left(f^n\left(D_{i,w}\right)\right)\leq A_n.$$
	
	Let $\pi_x:T_xM\to E_x^u$ denote the projection along $E_x^{cs}$. For $x\in S_{i,w}$, put $\mathcal{V}=T_xS_{i,w}$. Since the splitting is invariant,
	$$\pi_{f^n(x)}\circ Df_x^n|_\mathcal{V}=\left(Df_x^n|_{E_x^u}\right)\circ\pi_x|_\mathcal{V}.$$
	The projections $\pi_x$ are uniformly bounded, while the restrictions $\pi_x|_\mathcal{V}$ are uniformly invertible. Taking $u$-dimensional ($u:=\dim E^{u}$) Jacobians  therefore yields a constant $C>0$, independent of $n,i,w,x$, such that
	$$J_n^u\left(x\right)\leq C\operatorname{Jac}\left(Df_x^n|_\mathcal{V}\right).$$ As a consequence, Fubini's theorem and the change-of-variables formula give
	$$\int_M J_n^u\left(x\right)\,dx\leq C\sum_i\int\left(\int_{S_{i,w}}\operatorname{Jac}\left(Df_x^n|_{T_xS_{i,w}}\right)\,d\operatorname{Vol}_{S_{i,w}}\left(x\right)\right)dw=C\sum_i\int\operatorname{Vol}\left(f^n\left(S_{i,w}\right)\right)\,dw\leq CN A_n.$$
	
Taking logarithms, dividing by $n$, passing to the upper limit, and applying Proposition \ref{thm:cone-volume-growth}, we obtain the desired inequality.
\end{proof}

\vskip 5pt

Next we provide the proof of Theorem \ref{thm:main}.

\begin{proof}[Proof of Theorem \ref{thm:main}]
	
	The inequality in Theorem \ref{thm:main} is a combination of Proposition \ref{unstable entropy equal volume growth}, Proposition \ref{sup and lim order} and Lemma \ref{unstable entropy geq integrated volume growth}. 
	
	Let $\mu$ be a measure of maximal unstable entropy whose Lyapunov exponents along $E^{cs}$ are non-positive. Taking an ergodic component if necessary, we may assume that $\mu$ is ergodic. Then
	$$h_{\rm top}^{u}\left(f\right)=h_{\mu}^{u}\left(f\right)\leq h_{\mu}\left(f\right).$$
	By Saghin's entropy estimate \cite[Theorem 1.1]{Sag14},
	$$h_{\mu}\left(f\right)\leq\liminf_{n\to+\infty}\footnote{Saghin stated the inequality with a $\limsup$, but his proof actually yields the same inequality with $\liminf$.}\frac{1}{n}\log\int_{M}\left|\det\left(Df_x^n|_{E_x^u}\right)\right|\,dx+\dim\left(E^{cs}\right)\max\left\{\lambda_{\max}\left(E^{cs},\mu\right),0\right\},$$ where $\lambda_{\max}\left(E^{cs},\mu\right)$ denotes the maximal Lyapunov exponent of $\mu$ on $E^{cs}$.
	Since all Lyapunov exponents of $\mu$ along $E^{cs}$ are non-positive, the last term vanishes. Hence
	$$h_{\rm top}^{u}\left(f\right)\leq\liminf_{n\to+\infty}\frac{1}{n}\log\int_{M}\left|\det\left(Df_x^n|_{E_x^u}\right)\right|\,dx.$$
	Combining this with the inequality in the theorem yields the equality.
\end{proof}

\subsection{An example}
\begin{Example}\label{example >}
	There is a $C^{\infty}$ diffeomorphism on a compact manifold $M$ admitting a partially hyperbolic splitting $TM=E^{cs}\oplus E^{u}$ such that $$h^{u}_{\rm top}\left(f\right)>\lim_{n\to+\infty}\frac{1}{n}\log\int_{M}\left|\det\left(Df^{n}_{x}|_{E^{u}}\right)\right|\,dx.$$
\end{Example}
\begin{proof}
	
	Let $N$ be the unit tangent bundle of a closed surface of constant curvature $-1$, and let $\varphi^t:N\to N$ be its geodesic flow. We choose a Riemannian metric such that
	$$\left\|D\varphi^t|_{E_\varphi^u}\right\|=e^t,\qquad \left\|D\varphi^t|_{E_\varphi^s}\right\|=e^{-t}\qquad\left(t\geq0\right).$$
	Identify $\mathbb S^1$ with $\mathbb R/\mathbb Z$ and define
	$$g\left(\theta\right):=\theta+\frac{1}{4\pi}\sin\left(2\pi\theta\right)\pmod1.$$
	Then
	$$g'\left(\theta\right)=1+\frac12\cos\left(2\pi\theta\right),\qquad \frac12\leq g'\left(\theta\right)\leq\frac32,$$
	and $g\left(0\right)=0$, $g'\left(0\right)=3/2$. Set
	$$\tau\left(\theta\right):=\log\bigl(10g'\left(\theta\right)\bigr)$$
	and define the $C^\infty$ diffeomorphism
	$$f:N\times\mathbb S^1\to N\times\mathbb S^1,\qquad f\left(x,\theta\right):=\bigl(\varphi^{\tau\left(\theta\right)}\left(x\right),g\left(\theta\right)\bigr).$$
	Let $M:=N\times\mathbb S^1$.
	
	Let $X$ denote the vector field generating $\varphi^t$. The splitting
	$$E^u:=E_\varphi^u,\qquad E^{cs}:=\left(E_\varphi^s\oplus\mathbb RX\right)\oplus T\mathbb S^1$$
	is $Df$-invariant. Write
	$$T_n\left(\theta\right):=\sum_{j=0}^{n-1}\tau\left(g^j\theta\right).$$
	Then
	$$f^n\left(x,\theta\right)=\bigl(\varphi^{T_n\left(\theta\right)}\left(x\right),g^n\left(\theta\right)\bigr).$$
	Since
	$$T_n\left(\theta\right)=n\log10+\sum_{j=0}^{n-1}\log g'\left(g^j\theta\right)=n\log10+\log\left(g^n\right)'\left(\theta\right),$$
	we obtain
	$$\left|\det\left(Df^n_{\left(x,\theta\right)}|_{E^u}\right)\right|=e^{T_n\left(\theta\right)}=10^n\left(g^n\right)'\left(\theta\right).$$
	In particular,
	$$m\left(Df^n|_{E^u}\right)\geq5^n.$$
	
	Moreover,
	$$\left|T_n'\left(\theta\right)\right|\leq\left\|\tau'\right\|_\infty\sum_{j=0}^{n-1}\left|\left(g^j\right)'\left(\theta\right)\right|\leq C\left(\frac32\right)^n.$$
	For $v_s+aX+b\partial_\theta\in E^{cs}$, we have
	$$Df^n\left(v_s+aX+b\partial_\theta\right)=D\varphi^{T_n\left(\theta\right)}v_s+\bigl(a+bT_n'\left(\theta\right)\bigr)X+b\left(g^n\right)'\left(\theta\right)\partial_\theta.$$
	Hence
	$$\left\|Df^n|_{E^{cs}}\right\|\leq C\left(\frac32\right)^n,$$
	and therefore
	$$\frac{\left\|Df^n|_{E^{cs}}\right\|}{m\left(Df^n|_{E^u}\right)}\leq C\left(\frac3{10}\right)^n.$$
	Thus $f$ is partially hyperbolic with uniformly expanding bundle $E^u$.
	
	The unstable leaves of $f$ are precisely $W_\varphi^u\left(x\right)\times\left\{\theta\right\}$. Since the unstable Jacobian is constant along each such leaf,
	$$\operatorname{Vol}\left(f^n\left(W_\delta^u\left(x,\theta\right)\right)\right)=10^n\left(g^n\right)'\left(\theta\right)\operatorname{Vol}\left(W_\delta^u\left(x,\theta\right)\right).$$
	Since $\left(g^n\right)'\left(\theta\right)\leq\left(3/2\right)^n$ for every $\theta$, while
	$$\left(g^n\right)'\left(0\right)=\left(\frac32\right)^n,$$
	we obtain
	$$h_{\rm top}^u\left(f\right)=\lim_{n\to+\infty}\frac1n\log\sup_{z\in M}\operatorname{Vol}\left(f^n\left(W_\delta^u\left(z\right)\right)\right)=\log15.$$
	
	On the other hand, with respect to the normalized product Riemannian volume on $M$,
	$$\int_M\left|\det\left(Df^n_x|_{E^u}\right)\right|\,dx=10^n\int_{\mathbb S^1}\left(g^n\right)'\left(\theta\right)\,d\theta=10^n,$$
	since $g^n$ is an orientation-preserving circle diffeomorphism of degree one. Therefore
	$$\lim_{n\to+\infty}\frac1n\log\int_M\left|\det\left(Df^n_x|_{E^u}\right)\right|\,dx=\log10.$$
	Consequently,
	$$h_{\rm top}^u\left(f\right)=\log15>\log10=\lim_{n\to+\infty}\frac1n\log\int_M\left|\det\left(Df^n_x|_{E^u}\right)\right|\,dx.$$
\end{proof}

\begin{Remark*}
	$f$ admits a partially hyperbolic splitting $TM=E^s\oplus E^c\oplus E^u$, where $E^s:=E_\varphi^s$ and $E^c:=\mathbb RX\oplus T\mathbb S^1$, so that $\dim E^c=2$. Indeed, the preceding derivative formula also gives
	$$\left\|Df^n|_{E^s}\right\|\leq5^{-n},\qquad \left\|\left(Df^n|_{E^c}\right)^{-1}\right\|\leq C3^n,$$
	and hence
	$$\frac{\left\|Df^n|_{E^s}\right\|}{m\left(Df^n|_{E^c}\right)}\leq C\left(\frac35\right)^n.$$
\end{Remark*}

\section*{Appendix}

\begin{Lemma}\label{log-integral}
	Let $\left(X,\mathcal A,\mu\right)$ be a probability space and let $\varphi_n:X\to\left(0,+\infty\right)$ be measurable functions with $\left\{\frac1n\log\varphi_n\right\} $ being uniformly bounded.
	Assume
	$$\limsup_{n\to+\infty}\frac1n\int_X\log\varphi_n\,d\mu=\limsup_{n\to+\infty}\frac1n\log\int_X\varphi_n\,d\mu=h.$$
	Then there exists $n_k\to+\infty$ such that
	$$\frac1{n_k}\log\varphi_{n_k}\left(x\right)\longrightarrow h$$
	for $\mu$-almost every $x$.
\end{Lemma}

\begin{proof}
	Choose $n_k\to+\infty$ such that
	$$\frac1{n_k}\int_X\log\varphi_{n_k}\,d\mu\longrightarrow h.$$
	By Jensen's inequality and the assumption,
	$$\frac1{n_k}\log\int_X\varphi_{n_k}\,d\mu\longrightarrow h.$$
	For every $\varepsilon>0$,
	$$\mu\left\{x:\frac1{n_k}\log\varphi_{n_k}\left(x\right)\geq h+\varepsilon\right\}\leq\exp\left(\log\int_X\varphi_{n_k}\,d\mu-n_k\left(h+\varepsilon\right)\right)\longrightarrow0.$$
	By uniform boundedness of $\left\{\frac1n\log\varphi_n\right\} $, it follows that
	$$\int_X\left(\frac1{n_k}\log\varphi_{n_k}-h\right)^+d\mu\longrightarrow0.$$
	Together with
	$$\frac1{n_k}\int_X\log\varphi_{n_k}\,d\mu\longrightarrow h,$$
	we also obtain
	$$\int_X\left(h-\frac1{n_k}\log\varphi_{n_k}\right)^+d\mu\longrightarrow0.$$
	Hence
	$$\frac1{n_k}\log\varphi_{n_k}\longrightarrow h$$
	in $L^1\left(\mu\right)$. Passing to a subsequence gives the desired almost everywhere convergence.
\end{proof}

\section*{Acknowledgements}

ChatGPT was used for language-related tasks and for helping correct some reasoning errors. All mathematical ideas, proof strategies, and substantive arguments were developed and completed by the author.

Yuntao Zang is supported by National Key R\&D Program of China (2022YFA1005802) and NSFC (12201445, 12471186).

\vskip 5pt

\flushleft{\bf Yuntao Zang} \\
\small School of Mathematical Sciences,  Soochow University, Suzhou, 215006, P.R. China\\
\textit{E-mail:} \texttt{ytzang@suda.edu.cn}\\

\end{document}